\documentclass[a4paper,12pt]{amsart}
\usepackage{amssymb}
\usepackage{ifthen}
\nonstopmode \numberwithin{equation}{section}
\newtheorem{thm}[equation]{Theorem}
\newtheorem{cor}[equation]{Corollary}
\newtheorem{lem}[equation]{Lemma}
\newtheorem{prop}[equation]{Proposition}

\newtheorem{conj}{Conjecture}

\theoremstyle{definition}

\newtheorem{prob}[equation]{Problem}

\newenvironment{rem}{%
\bigskip
\noindent \textsl{{\sl Remark. }}}{\bigskip}
\newenvironment{rems}{%
\bigskip
\noindent \textsl{{\sl Remarks. }}}{\bigskip}

\newcounter {own}
\def\theown {\thesection       .\arabic{own}}

\newenvironment{pf}[1][]{%
 \vskip 3mm
 \noindent
 \ifthenelse{\equal{#1}{}}%
  {{\slshape Proof. }}%
  {{\slshape #1.} }%
 }%
{\qed
\medskip
}

\newcounter{alphabet}
\newcounter{tmp}

\newcommand{\IR}{{\mathbb R}}
\newcommand{\ID}{{\mathbb D}}

\newcommand{\IC}{{\mathbb C}}

\newcommand{\CC}{{\mathcal{C}}}

\def\be{\begin{equation}}
\def\ee{\end{equation}}

\newcommand{\bee}{\begin{enumerate}}
\newcommand{\eee}{\end{enumerate}}

\newcommand{\blem}{\begin{lem}}
\newcommand{\elem}{\end{lem}}
\newcommand{\bthm}{\begin{thm}}
\newcommand{\ethm}{\end{thm}}
\newcommand{\bcor}{\begin{cor}}
\newcommand{\ecor}{\end{cor}}
\newcommand{\beg}{\begin{examp}}
\newcommand{\eeg}{\end{examp}}
\newcommand{\begs}{\begin{examples}}
\newcommand{\eegs}{\end{examples}}
\newcommand{\bdefe}{\begin{defin}}
\newcommand{\edefe}{\end{defin}}
\newcommand{\bprob}{\begin{prob}}
\newcommand{\eprob}{\end{prob}}
\newcommand{\bei}{\begin{itemize}}
\newcommand{\eei}{\end{itemize}}

\newcommand{\bcon}{\begin{conj}}
\newcommand{\econ}{\end{conj}}
\newcommand{\bcons}{\begin{conjs}}
\newcommand{\econs}{\end{conjs}}
\newcommand{\bprop}{\begin{prop}}
\newcommand{\eprop}{\end{prop}}
\newcommand{\br}{\begin{rem}}
\newcommand{\er}{\end{rem}}
\newcommand{\brs}{\begin{rems}}
\newcommand{\ers}{\end{rems}}
\newcommand{\bo}{\begin{obser}}
\newcommand{\eo}{\end{obser}}
\newcommand{\bos}{\begin{obsers}}
\newcommand{\eos}{\end{obsers}}
\newcommand{\bpf}{\begin{pf}}
\newcommand{\epf}{\end{pf}}
\newcommand{\ba}{\begin{array}}
\newcommand{\ea}{\end{array}}
\newcommand{\beq}{\begin{eqnarray}}
\newcommand{\beqq}{\begin{eqnarray*}}
\newcommand{\eeq}{\end{eqnarray}}
\newcommand{\eeqq}{\end{eqnarray*}}

\newcommand{\ra}{\rightarrow}

\begin{document}
\bibliographystyle{amsplain}
\title{Regions of variability for a class of analytic and locally univalent functions defined by subordination}
\author{Bappaditya Bhowmik}
\address{Bappaditya Bhowmik,  Department of Mathematics,
Indian Institute of Technology Kharagpur, Kharagpur-721302, West Bengal, India.}
\email{bappaditya @ maths.iitkgp.ernet.in}

\date{ %\today
August 09, 2014; File: Revised-Bhow-Janowski.tex}
\begin{abstract}
In this article we consider a family $\mathcal{C}(A, B)$ of analytic and locally univalent functions on the open unit disc
$\ID=\{z :|z|<1\}$ in the complex plane that properly
contains the well-known Janowski class of convex univalent functions. In this article, we determine
the exact set of variability of $\log(f'(z_0))$ with fixed $z_0 \in \ID$ and $f''(0)$ whenever
$f$ varies over the class $\mathcal{C}(A, B)$.
\end{abstract}
\thanks{The author of this article would like to thank
ISIRD, SRIC,IIT Kharagpur \\(Ref.No.- IIT/SRIC/MATH/GNP/2014-15/61) for its financial support.}
%\thanks{}

\maketitle
\pagestyle{myheadings}
\markboth{B. Bhowmik  %and K.-J. Wirths
}{Regions of variability}

\section{Introduction and Preliminary Results}
Let $\IC$ be the complex plane. We use the following notations for open and closed discs with center $c$ and radius $r$ in the complex plane:
$\ID(c,r) := \{ z \in \mathbb{C} : |z-c| < r \}$ and $\overline{\ID}(c,r) := \{ z \in \mathbb{C} : |z-c| \leq  r \}$. We also denote
$\ID:= \ID(0,1)$ and $\overline{\ID} := \overline{\ID}(0,1)$. Let $\mathcal{A}$ be
the class of functions $f$ that are analytic in $\ID$ having the normalization $f(0)=0=f'(0)-1$
and $\mathcal{S}:=\{f\in \mathcal{A}\, :\, f~ \mbox{ is univalent in}~ \ID\}$. A function $f$ is said to be
{\it locally univalent} in $\ID$ if for any $z_0 \in \ID$, it is univalent in some neighborhood of $z_0$. A necessary and sufficient condition for
an analytic function to be locally univalent in $\ID$ is $f'(z) \not= 0$ in $\ID$.
We now state the definition of {\it subordination} which we need for our further
discussion. Let $h$ and $g$ be analytic functions in $\ID$. We call $h$ is {\it subordinate}
to $g$ (abbreviated symbolically as $h\prec g$) if there exists an analytic function $\xi: \ID\ra\ID$ such that
$\xi(0)=0$ and $h(z)=g(\xi(z))$. In special, if $g$ is one-one in $\ID$, then $h\prec g$ if and only if  $h(\ID)\subseteq g(\ID)$ with
$h(0)=g(0)$. We call the afore mentioned function $\xi : \ID \ra \ID$, which is analytic in $\ID$ with $\xi(0)=0$ as {\it Schwarz function} in the literature.

In 1973, Witold Janowski considered the following class of functions
\beqq
J_{A,B}:&=& \left\{f\in \mathcal{A}: \mbox{ there exists a Schwarz function}~\psi~ \mbox{in}~ \ID ~  \right.\\
&& \mbox{such that}~ 1+\frac{zf''(z)}{f'(z)}=\left. \frac{1+A\psi(z)}{1+B\psi(z)}\right\},
\eeqq
where $A$ and $B$ are real constants with $-1\leq A < B\leq 1$. We now observe that for $f\in\mathcal{A}$ and $-1\leq A < B<1$, $f\in J_{A, B}$   if and only if
the quantity $\left(1+\frac{zf''(z)}{f'(z)}\right)$ belongs to the disk which has the line segment $\left[\frac{1+A}{1+B}, \frac{1-A}{1-B}\right]$ as
a diameter, i.e.,
$$
1+\frac{zf''(z)}{f'(z)}\in \ID\left(\frac{1-AB}{1-B^2}, \frac{B-A}{1-B^2}\right)\, \forall\, z\in \ID.
$$
Since the above disk is contained in the right half plane, each $f\in J_{A, B}$ is convex univalent in $\ID$.
In case of $-1\leq A < B=1$, we have ${\rm Re}\, \left(1+\frac{zf''(z)}{f'(z)}\right)>0$ for $z\in \ID$, which again
forces each $f\in J_{A, B}$ to be convex univalent in $\ID$.
In particular, we have
$J_{A,B} \subsetneq \mathcal{S}$.
 We refer to the articles
\cite{Jan, Silverman} for details and many more interesting results about the Janowski class.
For a function $f\in J_{A, B}$, there exists a Schwarz function $\psi$ such that the following holds:
\beqq
&&1+\frac{zf''(z)}{f'(z)} = \frac{1+A\psi(z)}{1+B\psi(z)}\\
&\Rightarrow& 1+\frac{zf''(z)}{f'(z)}\prec  \frac{1+Az}{1+Bz}\\
&\Rightarrow& \frac{zf''(z)}{f'(z)}\prec \frac{(A-B)z}{1+Bz}.
\eeqq
Now letting $f'(z)=p(z)$, the last implication yields
\be\label{jan-eq-2}
\frac{zp'(z)}{p(z)}\prec \frac{(A-B)z}{1+Bz}=:\phi(z).
\ee
We now apply the result stated in \cite[Corolarry 3.1d.1, p. 76]{MM} to (\ref{jan-eq-2}) and get
$$
f'(z)=p(z) \,\prec\, \exp \int_{0}^{z}\frac{\phi(t)}{t}\, dt =: q_{A,B}(z).
$$
A little computation reveals that
\beqq
q_{A,B}(z) = \left\{
\ba{lll}
(1+Bz)^{\frac{A}{B}-1} \quad {\mbox{\rm for}}\quad B\neq 0, \\
\exp (Az)\quad {\mbox{\rm for}}\quad B=0.
\ea
\right.
\eeqq
The above discussion motivates us to consider the family $\mathcal{C}(A, B)$ as defined below:
%$$
%\mathcal{C}(A,B)= \left\{f\in\mathcal{S} : \log f' \prec \left(\frac{A}{B}-1\right)\log(1+Bz)\right\};
%$$
\beqq
\mathcal{C}(A, B)&=& \left\{f: f \mbox{ is analytic and locally univalent in } \ID\,
\mbox{ with} f(0)=f'(0)-1=0 \right.\\
&&\left.\mbox{ satisfying}\,
\log f'(z) \prec \left(\frac{A}{B}-1\right)\log(1+Bz), \, z\in \ID\right\};
\eeqq
where $-1\leq A < B\leq 1,~ B\neq 0$.
Here we clarify that, since $f'$ is a non-vanishing function on $\ID$, then $(f')^{\lambda}$ is well-defined and holomorphic once we
choose the determination $(f'(0))^{\lambda}=1$ for all $\lambda \in \IR$. In other words, we chose that branch of logarithm for which
$\log(f'(0))=0$. By definition we have $J_{A,B}\subseteq \mathcal{C}({A,B})$.
At this point we only confirm
that the class $\mathcal{C}({A,B})$ contains functions from $\mathcal{S}$.
We can see this if we consider
$A=0$ and $0<B\leq 1$. We then have from the definition of the class $\mathcal{C}({A,B})$,
$$
f'(z)=\frac{1}{1+ B\omega(z)},\, \mbox{where}~ \omega ~\mbox{is a Schwarz function}.
$$
A little computation shows that
$$
{\rm Re}\,f'(z)={\rm Re}\,\left(\frac{1}{1+ B\omega(z)}\right)>\frac{1}{2}~~ \mbox{ if}~ 0<B\leq 1.
$$
%where $\alpha (\simeq 0.543)$ is the real root of the equation $x^3 +x^2+ x-1=0$.
Therefore, an application of  Noshiro-Warschawski theorem yields that
$f$ is univalent in $\ID$. Hence, if $A=0$ and $0<B\leq 1$, we see that $\mathcal{C}({A,B})\subseteq \mathcal{S}$.
We remark here that the class $\mathcal{C}({A,B})$ defined above may or may not be a subclass of $\mathcal{S}$ for all $-1\leq A < B\leq 1,~ B\neq 0$.
We leave this problem open.

We next claim that $J_{A,B}\subsetneq \mathcal{C}(A, B)$. In order to establish our claim, we let $f\in \mathcal{C}(A, B)$.
Hence there exists a Schwarz function $\psi$ such that
$$
\log{f'(z)}=\left(\frac{A}{B}-1\right)\log(1+B\psi(z)),\quad z\in \ID.
$$
We deduce the following after differentiating the above expression,
$$
1+\frac{zf''(z)}{f'(z)}= 1+ \frac{(A-B)z\psi'(z)}{1+B\psi(z)}, \quad z\in \ID.
$$
In particular, letting $\psi(z)=z^2$, the above expression simplifies into
$$
1+\frac{zf''(z)}{f'(z)}=\frac{1+(2A-B)z^2}{1+Bz^2}, \quad z\in \ID.
$$
%Now the claim is true if we let $A=-1$ and $B=1$ and compare with the definition of the class $J_{-1,1}$.
The function $\frac{1+(2A-B)z^2}{1+Bz^2}$ maps $\ID$ onto the disc that have the line segment $\left[\frac{1+2A-B}{1+B}, \frac{1-2A+B}{1-B}\right]$
as a diameter. Since $\frac{1+2A-B}{1+B}< \frac{1+A}{1+B}$, it follows that $f \not\in J_{A, B}$. Hence, in particular $\mathcal{C}(A, B)\setminus J_{A, B} \neq \phi$.

Finally, for a function $f\in \mathcal{C}(A, B)$, we have
$$
f'(z) \prec (1+Bz)^{\frac{A}{B}-1}.
$$
Hence there exists an analytic function $\omega$, bounded by unity with $\omega(0)=0$ such that
$$
f'(z)=(1+B\omega(z))^{\frac{A}{B}-1};
$$
from which we have
\be\label{jan-1}
\omega(z)=\frac{(f'(z))^\frac{B}{A-B}-1}{B}.
\ee
Thereafter, we calculate
$$ \omega'(0)=\frac{f''(0)}{A-B}.
$$
Denoting $\omega'(0)=\lambda$, we have $|\lambda|\leq 1$ by virtue of the Schwarz Lemma.  Hence we fix the
second Taylor coefficients for functions in $\mathcal{C}(A, B)$ as
$$
f''(0)=\lambda (A-B);\quad \lambda\in\overline\ID.
$$
We now define the following class
\beq\nonumber
\CC_{\lambda}(A,B) &=& \{f\in \mathcal{C}(A, B) : f''(0)= \lambda(A-B) \},
~\mbox{ and set}\\\nonumber
V_{\lambda}(z_0,A,B)&=& \{\log f'(z_0): f\in \CC_{\lambda}(A,B)\};
\eeq
where $\lambda\in\overline\ID$ and $z_0\in\ID$ is an arbitrary but fixed complex number.
In this paper, we wish to determine the regions of variability $V_{\lambda}(z_0,A,B)$ of $\log f'(z_0)$ when $f$
ranges over the class $\CC_{\lambda}(A,B)$. Here we would like to mention that the similar problems on regions of variability
for some subclasses of $\mathcal{S}$ have already been studied in the articles \cite{BP,Yanagihara1,Yanagihara2} and some references
therein.

\section{Discussion on the set $V_{\lambda}(z_0, A, B)$ and
the Main result}\label{p3sec1}

We start this section by listing some basic properties of the set $V_{\lambda}(z_0, A, B)$:
\bee
\item We first note that $\CC_{\lambda}(A,B)$ is a compact subset of
the class of analytic functions in
$\ID$ endowed with the topology of uniform convergence on
compact subsets of $\ID$. As we know that for fixed $z_0\in \ID$, the map $\mathcal{C}_\lambda(A,B) \ni f\mapsto \log(f'(z_0))$ is continuous, therefore
we conclude that the set
$V_{\lambda}(z_0,A, B)$ is a compact subset of $\IC$.
%We remark here that $\CC_{\lambda}(A,B)$ is {\it not} a convex set.
 %and consequently $V_{\lambda}(z_0,A, B)$ {\it not} a convex subset of $\IC$.

\item If $|\lambda|=1$, then by applying Schwarz Lemma, we have $\omega(z)=\lambda z$ and consequently we deduce from (\ref{jan-1}) that
$$V_{\lambda}(z_0, A, B)=\left\{\left(\frac{A-B}{B}\right)\log(1+B\lambda z_0)\right\}.
$$
If $z_0=0$, then $V_{\lambda}(z_0, A, B)=\{0\}$.
Now, for $\lambda\in\ID$ and $a\in\overline\ID$, we introduce the following functions
\beq
\nonumber
\delta(z,\lambda) & = &\frac{z+\lambda}{1+\overline{\lambda}z}, \quad z\in\ID,\quad \mbox{and}\\
\label{p3eq4}\hspace{1.5cm}
F_{a,\lambda}(z) & = & \int_0^z \left(1+B\zeta \delta(a\zeta,\lambda)\right)^{\frac{A-B}{B}} d\zeta.
\eeq
The above integral equation yields
$$
\log F'_{a,\lambda}(z)= \left({\frac{A-B}{B}}\right)\log\left(1+Bz \delta(az,\lambda)\right).
$$
Thus we have
$$
\log F'_{a,\lambda}(z) \prec \left(\frac{A}{B}-1\right)\log(1+Bz)
$$
and $F''_{a,\lambda}(0)=\lambda(A-B)$, showing that $F_{a,\lambda}\in \CC_{\lambda}(A, B)$. We observe that,
for a fixed $\lambda\in\ID$ and $z_0\in\ID\setminus\{0\}$, the function $a \mapsto \log F'_{a,\lambda}(z_0)$ is a non-constant
analytic function and hence is an open mapping. From this, we see that
$$\log F'_{0,\lambda}(z_0)= \left(\frac{A-B}{B}\right)\log(1+B\lambda z_0)
$$
is an interior point of $\{\log F'_{a,\lambda}(z_0):\,a\in\ID\}\subset V_{\lambda}(z_0, A, B)$.
%Since $V_{\lambda}(z_0, A, B)$ is a compact subset of $\IC$
%and also has nonempty interior, the boundary $\partial{V_{\lambda}(z_0, A, B)}$
%is a Jordan curve and $V_{\lambda}(z_0,A, B)$ is the union of
%$\partial{V_{\lambda}(z_0,A, B)}$ and the domain enclosed by this boundary.

\item We claim that $V_{\lambda}(e^{i\theta}z_0, A, B)= V_{\lambda e^{i\theta}}(z_0, A, B)$ for $\theta \in \IR$. This will be established
if we can show that $f\in \CC_{\lambda}(A,B)$   if and only if $e^{-i\theta}f(e^{i\theta}z)\in \CC_{\lambda e^{i\theta}}(A,B)$.
To prove this, let $f\in \CC_{\lambda}(A,B)$ and define
$$g(z)= e^{-i\theta}f(e^{i\theta}z),\quad z\in \ID.
$$
A little computation reveals that
\beqq
g''(0) &=& e^{i\theta}f''(0)=(A-B)\lambda e^{i\theta}\quad \mbox{ and}\\
\log g'(z)&=& \log f'(e^{i\theta}z)= \left(\frac{A-B}{B}\right)\log(1+B\omega({e^{i\theta}z)}).
\eeqq
Hence we have
$$\log g'(z) \prec \left(\frac{A-B}{B}\right)\log(1+Bz),
$$
proving $g(z)\in \CC_{\lambda e^{i\theta}}(A,B)$.
Conversely,  let
$$
g(z)= e^{-i\theta}f(e^{i\theta}z)\in \CC_{\lambda e^{i\theta}}(A,B).
$$
Hence we have
$$g''(0)= \lambda e^{i\theta} (A-B)= e^{i\theta}f''(0),
$$
resulting $f''(0)= \lambda(A-B)$. Next we observe the following chain of implications:
\beqq &&\log g'(z) \prec \left(\frac{A-B}{B}\right)\log(1+Bz)\\
&&\Rightarrow \log f'(e^{i\theta}z) \prec \left(\frac{A-B}{B}\right)\log(1+Bz)\\
&&\Rightarrow \log f'(z) \prec \left(\frac{A-B}{B}\right)\log(1+Bz).
\eeqq
This proves $f\in \CC_\lambda(A,B)$.
\eee
In view of the afore mentioned  properties of the set $V_{\lambda}(z_0,A,B)$, it is sufficient to determine $V_{\lambda}(z_0,A,B)$
for $\lambda \in [0,1)$, as the case $|\lambda| =1$ is completely described in the Item~$(2)$ above.  We now state our main theorem of the
paper.
\bthm\label{p3th1}
For $\lambda\in [0,1)$ and $z_0\in\ID\setminus \{0\}$, we have
$$
V_{\lambda}(z_0, A, B)=\left\{\left(\frac{A-B}{B}\right)\log \left(c(z_0,\lambda)+ar(z_0,\lambda)\right): |a|\leq 1\right\},
$$
where
\beqq
c(z_0,\lambda)
 & = & \frac{1-\lambda^2|z_0|^2+\lambda B(1-|z_0|^2)z_0}{1-\lambda^2|z_0|^2},\quad \mbox{and}\\
 r(z_0,\lambda) & = & \frac{|B|(1-\lambda^2)|z_0|^2}{1-\lambda^2|z_0|^2}.
\eeqq
The boundary $\partial{V_{\lambda}(z_0, A, B)}$ of the above set is the Jordan curve given
by
\beq\label{p3eq7}
(-\pi,\pi]\ni \theta &\mapsto& \log F'_{e^{i\theta},\lambda}(z_0)\\
 &=& \left(\frac{A-B}{B}\right) \log\left(1+Bz_0 \delta(e^{i\theta}z_0,\lambda)\right).\nonumber
\eeq
If  $\log f'(z_0)=\log F'_{e^{i\theta},\lambda}(z_0)$ for some
$f\in \CC_{\lambda}(A, B)$ and $\theta\in (-\pi,\pi]$, then $f(z)=
F_{e^{i\theta},\lambda}(z)$. Here $F_{e^{i\theta},\lambda}(z)$ is
given by $(\ref{p3eq4})$ with $a=e^{i\theta}$.
\ethm

We prove the Theorem~\ref{p3th1}  in the following Section.

\section{Description of the set $V_{\lambda}(z_0,A,B)$}
At the beginning, we establish the following result which will help us to achieve our main goal, i.e. to  prove the Theorem~\ref{p3th1}.
\bprop\label{p3prop1}
For $f\in \CC_{\lambda}(A,B)$ with $\lambda\in [0,1)$,
we have
\be\label{p3eq8}
\left|(f'(z))^{\frac{B}{A-B}}-c(z,\lambda)\right|\leq r(z,\lambda),
\quad z\in\ID,
\ee
where
\beqq
c(z,\lambda)
& =& \frac{1-\lambda^2|z|^2+\lambda B(1-|z|^2)z}{1-\lambda^2|z|^2}\quad \mbox{and}\\
 r(z,\lambda) & = & \frac{|B|(1-\lambda^2)|z|^2}{1-\lambda^2|z|^2}.
\eeqq
For each $z\in\ID\setminus\{0\}$, equality holds in $(\ref{p3eq8})$ if and
only if $f=F_{e^{i\theta},\lambda}$ for some
$\theta\in\IR$.
\eprop
\bpf
Let $f\in \CC_{\lambda}(A,B)$. Then, as stated in the introduction, there exists an $\omega$ such that
$$
\omega(z)=\frac{(f'(z))^\frac{B}{A-B}-1}{B}, \quad z\in \ID.
$$
Since
$\omega$ is a Schwarz function that satisfies the condition $\omega'(0)=\lambda$,
it follows from the Schwarz lemma that
\be\label{p3eq9}
\left|\frac{\frac{\omega(z)}{z}-\lambda}{1-\lambda\frac{\omega(z)}{z}}\right|
\leq|z|,\quad z\in\ID.
\ee
Now we insert the above expression of $\omega(z)$ in the inequality (\ref{p3eq9}) and we get
\be\label{p3eq10}
\left|\frac{(f'(z))^{\frac{B}{A-B}}-M(z,\lambda)}{\lambda(f'(z))^{\frac{B}{A-B}}-N(z,\lambda)}
\right| \leq |z|,\quad z\in\ID ,
\ee
where
\beq\label{p3eq11}
\left\{
\ba{lll}
M(z,\lambda)&=& 1+\lambda Bz\\[3mm]%\vspace*{.5cm}
N(z,\lambda)&=& \lambda+Bz.
\ea
\right.
\eeq
A  calculation shows that the inequality (\ref{p3eq10}) is
equivalent to the following inequality
\be\label{p3eq12}
\left|(f'(z))^{\frac{B}{A-B}}- \frac{M -\lambda|z|^2 N}{1-\lambda^2|z|^2}\right| \leq \frac{|\lambda M-N||z|}{1-\lambda^2|z|^2}, \quad z\in\ID.
\ee
We now simplify the following expressions using (\ref{p3eq11}) as
\beqq
%\label{p3eq13}
\frac{M -\lambda|z|^2 N}{1-\lambda^2|z|^2}
 & = & \frac{1-\lambda^2|z^2|+\lambda B(1-|z|^2)z}{1-\lambda^2|z|^2},~\mbox{ and}\\
%\label{p3eq14}
\frac{|\lambda M-N||z|}{1-\lambda^2|z|^2}
 & = & \frac{|B|(1-\lambda^2)|z|^2}{1-\lambda^2|z|^2}.
\eeqq
Therefore, the inequality (\ref {p3eq8}) follows from the last two
equalities and (\ref{p3eq12}). We now turn to the case of establishing sharpness of the inequality (\ref{p3eq8}).
Our first aim is to prove that the equality occurs for any $z\in\ID$ in (\ref{p3eq8}) whenever
$f=F_{e^{i\theta},\lambda}$, for some $\theta\in\IR$. To this end, we first show that if  $f=F_{e^{i\theta},\lambda}$ then
$$
\left({{F'}_{e^{i\theta},\lambda}}(z)\right)^{\frac{B}{A-B}}= 1+ B\left(\frac{e^{i\theta}z+\lambda}{1+\lambda ze^{i\theta}}\right)z,\quad \theta\in \IR.
$$
A straightforward computation reveals that
\beq\label{jan-eq-1}
&&\left({{F'}_{e^{i\theta},\lambda}}(z)\right)^{\frac{B}{A-B}}-c(z,\lambda) \\\nonumber
&=&1+ B\left(\frac{e^{i\theta}z+\lambda}{1+\lambda e^{i\theta}z}\right)z-
\frac{1-\lambda^2|z^2|+\lambda B(1-|z|^2)z}{1-\lambda^2|z|^2}\\ \nonumber
&=&\left(\frac{B(1-\lambda^2)z^2}{1-\lambda^2|z|^2}\right)~\left(\frac{e^{i\theta}+\lambda\overline z}{1+z\lambda e^{i\theta}}\right), \quad z\in \ID.
\eeq
The first part of the case of equality follows from (\ref{jan-eq-1}) by noting that
$$
\left|\frac{e^{i\theta}+\lambda\overline z}{1+z\lambda e^{i\theta}}\right|=1,\, z\in \ID;
$$
for the aforesaid range of values of $\theta$ and $\lambda$.
Conversely, if the equality holds for some
$z\in\ID\setminus\{0\}$ in (\ref {p3eq8}), then the
equality must hold in (\ref{p3eq9}). Consequently from the Schwarz
lemma, there exists $\theta\in \IR$ such that
$\omega(z)=z\delta(e^{i\theta}z,\lambda)$ for all
$z\in\ID$, i.e. we have
$$
\frac{(f'(z))^\frac{B}{A-B}-1}{B}=z\left(\frac{e^{i\theta}z+\lambda}{1+ \lambda e^{i\theta}z}\right),
$$
which after solving for $f$ results $f(z)=F_{e^{i\theta},\lambda}(z)$, $z\in \ID$. This completes the proof of the Proposition.
\epf

Next, we get the following useful estimate in the case of $\lambda =0$ that one may look for.
\bcor\label{p3cor1a}
Let $f\in \mathcal{C}_{0}( A, B)$. Then we have the following sharp inequality:
\be\label{p3eq19}
\left|(f'(z))^{\frac{B}{A-B}}-1\right|\leq |B||z|^2, \quad z\in \ID.
\ee
\ecor

We are now ready to deliver a proof of our main theorem after all these preparation.
\smallskip
\bpf[{\bf Proof of the Theorem~\ref{p3th1}}]
We begin the proof by defining the following class of functions:
$$
H^\infty ( \mathbb{D} )
= \{ w : w \text{ is analytic in $\mathbb{D}$ with $ | w(z)| \leq 1 $, $z\in\mathbb{D}$ } \} .
$$
We observe that $f \in \mathcal{C} (A, B)$ if and only if there exists $\varphi \in H^\infty ( \mathbb{D})$ such that
$$
f'(z) = (1+Bz\varphi (z))^{\frac{A-B}{B}}, z\in \ID.
$$
Furthermore, $f \in \mathcal{C}_\lambda (A, B)$
if and only if
there exists $w \in H^\infty ( \mathbb{D})$
such that
\beqq
f'(z) &=&\left( 1+Bz \frac{zw(z)+\lambda }{1+ \overline{\lambda}z w(z) }\right)^{\frac{A-B}{B}} \\
&=&\left( 1+Bz \delta (z w(z),\lambda )\right)^{\frac{A-B}{B}},\, z\in \ID.
\eeqq
%Now we consider the following class of functions:
%\begin{align*}
%  \mathcal{G}_\lambda (A,B)
%  =&
%  \left\{ f : \text{$f$ is analytic in $\mathbb{D}$
% and there exists $\omega \in H_1^\infty ( \mathbb{D} )$
% such that  } \right. \\
%  &
%   \left.
% f'(z) =  \left( 1+Bz \delta (z\omega(z),\lambda )
%  \right)^{\frac{A-B}{B} }, z\in \mathbb{D} \right\} .
%\end{align*}
%We notice that $\mathcal{C}_\lambda (A, B) = \mathcal{G}_\lambda (A,B)$.
Therefore we get
\begin{align*}
 V_\lambda (z_0,A,B)
 =&
 \left\{ \log f'(z_0) : f \in \mathcal{C}_\lambda (A,B) \right\}
\\
 =&
% \left\{ \log f'(z_0) : f \in \mathcal{G}_\lambda (A,B) \right\}
%\\
% =&
 \left\{  \left(\frac{A-B}{B}\right)
 \log \left( 1+Bz_0 \delta (z_0 w(z_0),\lambda ) \right)
 : w\in H^\infty ( \mathbb{D} )
\right\}
\\
 =&
 \left(\frac{A-B}{B}\right)
 \log \left( 1+Bz_0 \delta (z_0 \overline{\mathbb{D}},\lambda ) \right) .
\end{align*}
Combining the above and the relation
$$
1+Bz_0 \delta (z_0 \overline{\mathbb{D}},\lambda )= \overline{\mathbb{D}}(c(z_0,\lambda ), r(z_0, \lambda ) ),
$$
which easily follows from the proof of the Proposition~\ref{p3prop1},
we obtain
$$
 V_\lambda (z_0,A,B) =
 \left\{
 \left(\frac{A-B}{B}\right)
 \log \left( c(z_0,\lambda ) + a r(z_0, \lambda ) \right) : a \in
 \overline{\mathbb{D}}\right\} .
$$
Since the function $\log (1+w)$ is convex univalent in $\ID$,
the set in right hand side of the above equation
is a convex closed Jordan domain and its boundary curve is given by
$$
(- \pi , \pi ] \ni \theta
\mapsto \left(\frac{A-B}{B}\right) \log \left( 1+Bz_0 \delta (z_0 e^{i \theta },\lambda ) \right)
= \log F_{e^{i \theta }, \lambda }'(z_0) .
$$
This completes the proof of the Theorem.
\epf

\bigskip
{\bf Acknowledgement:} The author thanks K-J. Wirths for his suggestions. The author would also
like to thank the referee for his careful reading of the paper and inputs in the proof of the Theorem~\ref{p3th1}.

\end{document}